\documentclass[a4paper,10pt,reqno, english]{amsart}  
\usepackage[utf8]{inputenc}
\usepackage[T1]{fontenc}
\usepackage{amsmath,amsthm}
\usepackage{amsfonts,amssymb,enumerate}
\usepackage{url,paralist}
\usepackage{mathtools}  
\usepackage[colorlinks=true,urlcolor=blue,linkcolor=red,citecolor=magenta]{hyperref}
\usepackage{enumerate}
\usepackage{anysize}
\theoremstyle{plain}
\newtheorem{theorem}{Theorem}[section]

\newtheorem{corollary}[theorem]{Corollary}

\newtheorem*{theorem*}{Theorem}

\newtheorem*{claim*}{Claim}

\theoremstyle{definition}

\newcommand{\R}{{\rm I\!R}}

\newcommand{\conv}{\operatorname{conv}}
\newcommand{\card}{\operatorname{card}}
\newcommand{\aff}{\operatorname{aff}}
\newcommand{\vset}{\operatorname{vert}}
\newcommand{\spann}{\operatorname{span}}

\begin{document}

\title[Plus minus analogues for affine Tverberg type results]
{Plus minus analogues for affine Tverberg type results}
 
\author[Blagojevi\'c]{Pavle V. M. Blagojevi\'{c}}
\address[Pavle V. M. Blagojevi\'{c}]{Institut f\" ur Mathematik, FU Berlin, Arnimallee 2, 14195 Berlin, Germany\hfill\break%
\mbox{\hspace{4mm}}Mathematical Institute SASA, Knez Mihailova 36, 11000 Beograd, Serbia}
\email{blagojevic@math.fu-berlin.de} 
\author[Ziegler]{G\"unter M. Ziegler}  
\address[G\"unter M. Ziegler]{Inst. Math., FU Berlin, Arnimallee 2, 14195 Berlin, Germany} 
\email{ziegler@math.fu-berlin.de}

\thanks{The research by Pavle V. M. Blagojevi\'{c} leading to these results has received funding from DFG via the Collaborative Research Center TRR~109 ``Discretization in Geometry and Dynamics,'' and the grant ON 174024 of the Serbian Ministry of Education and Science.\newline\indent
The research by G\"unter M. Ziegler leading to these results has funding from DFG via the Collaborative Research Center TRR~109 ``Discretization in Geometry and Dynamics.'' 
}

\date{}

\dedicatory{Dedicated to the memory of Branko Gr\"unbaum}

\maketitle


\begin{abstract}
The classical 1966 theorem of Tverberg with its numerous variations was and still is a motivating force behind many important developments in convex and computational geometry as well as a testing ground for methods from equivariant algebraic topology.
In 2018, B\'ar\'any and Sober\'on presented a new variation, the ``Tverberg plus minus theorem.''
In this paper, we give a new proof of the Tverberg plus minus theorem, by using a projective transformation.
The same tool allows us to derive plus minus analogues of all known affine Tverberg type results.
In particular, we prove a plus minus analogue of the optimal colored Tverberg theorem. 
 
\end{abstract}
 
\section{Introduction and the statement of main result}
\label{sec:Introduction}

A seminal 1966 result of Helge Tverberg \cite{Tverberg1966} states that for any integers $d\ge1$ and $r\ge2$, with $N=(r-1)(d+1)$, every affine map of the $N$-dimensional simplex $\Delta_N$ into the Euclidean space $\R^d$ identifies points from $r$ vertex-disjoint faces of $\Delta_N$.
Over the years, along with numerous applications, Tverberg's original theorem was extended in many different directions. 
For details on the history, relevance and connections to equivariant topology consult \cite{Matousek2008} \cite{Ziegler2011} \cite{Barany2016} and \cite{Barany2018s}; or for complete proofs of all relevant results see for example \cite{Blagojevic2017}. 
 
We now consider a recent variation of the classical Tverberg's theorem, the ``Tverberg plus minus theorem'' of Imre B\'ar\'any and Pablo Sober\'on \cite{Barany2018}.
In the following, for a face $\mu$ of the simplex $\Delta_N$ we set $\mu'$ to be the complementary face of the face $\mu$ in $\Delta_N$, that is,
\[
\vset(\mu)\cup\vset(\mu')=\vset(\Delta_N)
\qquad\text{and}\qquad
\vset(\mu)\cap\vset(\mu')=\emptyset.
\]
Here $\vset(\mu)$ denotes the set of vertices of the face $\mu$ of the simplex $\Delta_N$.

\begin{theorem}[The Tverberg plus minus theorem of B\'ar\'any and Sober\'on {\cite[Thm.\,1.3]{Barany2018}}]
\label{th:TverbergPlusMinus}
Let $d\ge1$ and $r\ge2$ be integers with $N=(r-1)(d+1)$, and let $a\colon \Delta_N\longrightarrow\R^d$ be an affine map.
Furthermore, let $\mu$ be a face of the simplex $\Delta_N$ of dimension at most $r-2$  with $a(\mu)\cap a(\mu')=\emptyset$.
Then there exist $r$ pairwise disjoint proper faces $\sigma_1,\dots,\sigma_r$ of $\Delta_N$ and there exists a point $b\in\R^d$ in the intersection  $\aff(a(\sigma_1))\cap\dots\cap \aff(a(\sigma_r))$ having a representation by  $r$ affine combinations
\[
b=\sum_{v\in \vset(\sigma_1)}\alpha_va(v)=\dots=\sum_{v\in \vset(\sigma_r)}\alpha_va(v),
\qquad
1=\sum_{v\in \vset(\sigma_1)}\alpha_v=\dots=\sum_{v\in \vset(\sigma_r)}\alpha_v,
\] 
such that 
\[
(v\in\vset(\mu)\  \Longrightarrow \  \alpha_v\leq 0)\quad\text{and}\quad (v\in\vset(\mu')\  \Longrightarrow \  \alpha_v\geq 0).
\]
\end{theorem}
 
Our statement of the Tverberg plus minus theorem differs from the original version in \cite[Thm.\,1.3]{Barany2018} because we have  dropped a general position assumption on the affine map $a\colon \Delta_N\longrightarrow\R^d$, that all $(N+1)d$ coordinates of the point set $\{ a(v) : v\in\vset(\Delta_N)\}$ are algebraically independent. 
With this extra condition the point $b$ becomes unique \cite[Prop.\,1.1]{Barany2018}, and the condition on the signs of the coefficients in the representations of $b$ becomes: $(v\in\vset(\mu)\ \Leftrightarrow \alpha_v < 0)$  and $(v\in\vset(\mu')\  \Leftrightarrow   \alpha_v> 0)$.
All the results of this paper will be stated and proven without any extra assumptions on the affine map $a\colon \Delta_N\longrightarrow\R^d$, but, as in the work of B\'ar\'any and Sober\'on, all the results can be sharpened in the same way if the algebraic independence of the coordinates of the set $\{ a(v) : v\in\vset(\Delta_N)\}$ is assumed.
 
Using an appropriate projective transformation, we show that the Tverberg plus minus theorem is a corollary of the classical Tverberg theorem.
The new proof of Theorem \ref{th:TverbergPlusMinus}, given in Section \ref{sec:ProofOfTverbergPlusMinus}, allows us to directly derive plus minus analogues of literally all known affine Tverberg type results.
On the other hand, it is not clear whether our proof of Theorem \ref{th:TverbergPlusMinus} implies plus minus analogues of integer Tverberg type results \cite{Loera2019}. 
Briefly, in our approach we begin by fixing parameters and assumptions depending on the Tverberg type result for which we want to derive a plus minus analogue.
For example, 
\begin{compactitem}[\quad---]
\item if we want to prove a plus minus analogue of a basic Tverberg type result, the parameters we fix are $d\geq 1$ and $r\geq 2$ integers, $N:=(r-1)(d+1)$, and $a\colon\Delta_N\longrightarrow\R^d$ an affine map;
\item if we want to prove a plus minus analogue of a optimal colored Tverberg type result, the parameters we fix are $d\geq 1$ an integer, $r\geq 2$ a prime, $N:=(r-1)(d+1)$, $a\colon\Delta_N\longrightarrow\R^d$ an affine map, and a coloring of the set of vertices of the simplex $\Delta_N$ by $m\geq 1$ colors such that each color class does no contain more that $r-1$ vertices;
\item if we want to prove a plus minus analogue of a generalized Van Kampen--Flores type result, the parameters we fix are $d\geq 1$ an integer, $r\geq 2$ a prime power, $N:=(r-1)(d+2)$, and $a\colon\Delta_N\longrightarrow\R^d$ an affine map.
\end{compactitem}
Then, using the assumption about the existence of the face $\mu$ of the simplex $\Delta_N$ with the properties that $a(\mu)\cap a(\mu')=\emptyset$ and $\dim(\mu)\leq r-2$, we 
\begin{compactitem}[\quad---]
\item transform the affine map $a\colon\Delta_N\longrightarrow\R^d\lhook\joinrel\longrightarrow\R^{d+1}$, via a permissible projective transformation of $\R^{d+1}$, into an affine map $\widehat{a}\colon \Delta_N\longrightarrow E'\lhook\joinrel\longrightarrow\R^{d+1}$ from the same simplex into a different Euclidean space of the same dimension,
\item to which we apply our chosen Tverberg type theorem.
This yields the corresponding, always affine, plus minus analogue result.
\end{compactitem}
In this way we get plus minus analogues of the optimal colored Tverberg theorem, the weak colored Tverberg theorem, the generalized Van Kampen--Flores theorem, the equal barycentric coordinates Tverberg theorem of Sober\'on, and the admissible-prescribable balanced theorem.
 
In particular, we get a plus minus analogue of the optimal colored Tverberg theorem of Blagojevi\'c, Matschke, and Ziegler \cite[Thm.\,2.1]{Blagojevic2009} \cite[Thm.\,2.1]{Blagojevic2011-2}.
The original result was obtained in the process of proving the B\'ar\'any-Larman colored Tverberg conjecture \cite{Barany1992} for the case when $r+1$ is a prime.

\begin{theorem}[Plus minus optimal colored Tverberg theorem]
\label{th:OptimalColoredTverbergPlusMinus}
Let $d\ge1$ be an integer, $r\ge2$ be a prime with $N=(r-1)(d+1)$, and let $a\colon \Delta_N\longrightarrow\R^d$ be an affine map.
Furthermore, let $\mu$ be a face of the simplex $\Delta_N$ of dimension at most $r-2$ with $a(\mu)\cap a(\mu')=\emptyset$, and let the vertex set of the simplex $\Delta_N$ be partitioned into $m+1$ parts, called ``color classes,''  
\[
\vset(\Delta_N)=C_0\sqcup\dots\sqcup C_m
\]
such that the cardinality of each color class does not exceed $r-1$, that is $\card(C_i)\leq r-1$ for all $0\leq i\leq m$.
Then there exist $r$ pairwise disjoint proper faces $\sigma_1,\dots,\sigma_r$ of $\Delta_N$, and there exists a point $b\in\R^d$ in the intersection  $\aff(a(\sigma_1))\cap\dots\cap \aff(a(\sigma_r))$ having a representation by  $r$ affine combinations
\[
b=\sum_{v\in \vset(\sigma_1)}\alpha_va(v)=\dots=\sum_{v\in \vset(\sigma_r)}\alpha_va(v),
\qquad
1=\sum_{v\in \vset(\sigma_1)}\alpha_v=\dots=\sum_{v\in \vset(\sigma_r)}\alpha_v,
\] 
such that  
\[
(v\in\vset(\mu)\  \Longrightarrow \  \alpha_v\leq 0)\quad\text{and}\quad (v\in\vset(\mu')\  \Longrightarrow \  \alpha_v\geq 0),
\]
and in addition
\[
\card (C_i\cap \vset(\sigma_j))\leq 1
\]
for every $0\leq i\leq m$ and every $1\leq j\leq r$. 
That is, all the faces $\sigma_1,\dots,\sigma_r$ are ``rainbow faces'' -- no two vertices of each $\sigma_i$ are colored with the same color.
\end{theorem}

\begin{figure}[ht]
\includegraphics{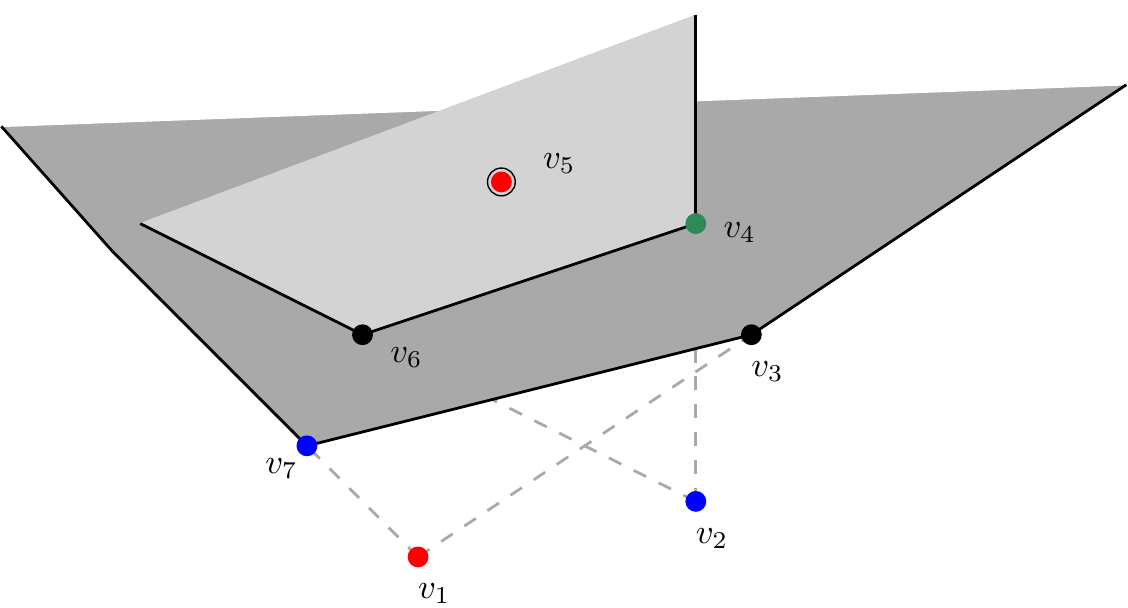}
\caption{\small An illustration of Theorem \ref{th:OptimalColoredTverbergPlusMinus} for the collection of colored points in the plane where $d=2$, $r=3$, $N=6$, $C_0=\{v_1,v_5\},C_1=\{v_2,v_7\},C_2=\{v_3,v_6\},C_3=\{v_4\}$, $\vset(\mu)=\{v_1,v_2\}$, $\vset(\mu')=\{v_3,v_4,v_5,v_6,v_7\}$, $\sigma_1=\{v_2,v_4,v_6\}$, $\sigma_2=\{v_1,v_3,v_7\}$, and $\sigma_3=\{v_5\}$.}	
\end{figure}

Identifying one color class with the set of vertices of the distinguished face $\mu$ of the simplex $\Delta_N$ yields the 
following strengthening of the original Tverberg plus minus theorem.

\begin{corollary}
Let $d\ge1$ be an integer, $r\ge2$ be a prime with $N=(r-1)(d+1)$, and let $a\colon \Delta_N\longrightarrow\R^d$ be an affine map.
Furthermore, let $\mu$ be a face of the simplex $\Delta_N$ of dimension at most $r-2$ with $a(\mu)\cap a(\mu')=\emptyset$.	
Then there exist $r$ pairwise disjoint proper faces $\sigma_1,\dots,\sigma_r$ of $\Delta_N$, and there exists a point $b\in\R^d$ in the intersection  $\aff(a(\sigma_1))\cap\dots\cap \aff(a(\sigma_r))$ having a representation by  $r$ affine combinations
\[
b=\sum_{v\in \vset(\sigma_1)}\alpha_va(v)=\dots=\sum_{v\in \vset(\sigma_r)}\alpha_va(v),
\qquad
1=\sum_{v\in \vset(\sigma_1)}\alpha_v=\dots=\sum_{v\in \vset(\sigma_r)}\alpha_v,
\] 
such that   
\[
(v\in\vset(\mu)\  \Longrightarrow \  \alpha_v\leq 0)\quad\text{and}\quad (v\in\vset(\mu')\  \Longrightarrow \  \alpha_v\geq 0),
\]
and in addition
\[
\card (\vset(\mu)\cap \vset(\sigma_j))\leq 1
\]
for every $1\leq j\leq r$.
That is, in each affine presentation $\sum_{v\in \vset(\sigma_i)}\alpha_va(v)$ of $b$ there exists at most one vertex $v\in\vset(\sigma_i)$ with the corresponding coefficient $\alpha_v$ being non-positive.
\end{corollary}  

In the next section we give in parallel proofs of Theorem \ref{th:TverbergPlusMinus} and Theorem \ref{th:OptimalColoredTverbergPlusMinus}.
The plus minus analogues of the weak colored Tverberg theorem, the generalized Van Kampen--Flores theorem, the equal barycentric coordinates Tverberg theorem of Soberón, and the  admissible-prescribable balanced Tverberg theorem, can be proved in the exact same way.  
 
\subsection*{Acknowledgement} 
We are grateful to Imre B\'ar\'any, Florian Frick, and Pablo Sober\'on for useful discussions and comments.
Furthermore, we thank the referees, Peter Landweber and Pablo Sober\'on for careful reading of the manuscript, valuable suggestions and corrections.

\section{Proofs of Theorem \ref{th:TverbergPlusMinus} and Theorem \ref{th:OptimalColoredTverbergPlusMinus}}
\label{sec:ProofOfTverbergPlusMinus}

In this section we simultaneously prove Theorem \ref{th:TverbergPlusMinus} and Theorem \ref{th:OptimalColoredTverbergPlusMinus}.
All the remaining analogues can be obtained in exactly the same way by modifying initial parameters and assumptions accordingly.

\medskip
Let $d\ge1$ be an integer, let $r\ge2$ be either an arbitrary integer or just a prime, and set $N=(r-1)(d+1)$.
Let $a\colon \Delta_N\longrightarrow \R^d$ be an affine map.
Denote by $V:=\vset(\Delta_N)=\{v_0,\dots,v_N\}$ the set of vertices of the simplex $\Delta_N$, and by $I:=\{0,\dots, N\}$ the corresponding index set.
In addition, when considering Theorem \ref{th:OptimalColoredTverbergPlusMinus}, we assume that the vertex set $V$ of the simplex $\Delta_N$ is partitioned into $m+1$ color classes $V=C_0\sqcup\dots\sqcup  C_m$ such that the cardinality of each color class does not exceed $r-1$.

\medskip
Let $\mu$ be a face of $\Delta_N$ with the property that $\dim(\mu)\leq r-2$ and $a(\mu)\cap a(\mu')=\emptyset$. 
Consequently, there exists an affine hyperplane $H$ in $\R^d$ that strictly separates $a(\mu)$ and $a(\mu')$.
Let $H$ be given by
\begin{equation*} 
	H:=\{ x\in \R^d : \langle x,w\rangle =\alpha \}
\end{equation*}
where the vector $w\in\R^d{\setminus}\{0\}$ and the scalar $\alpha\in\R$ are chosen in such a way that 
\begin{equation}
	\label{separation}
	a(\mu)\subseteq \mathrm{int}(H^-)=\{ x\in \R^d : \langle x,w\rangle <\alpha \} 
	\qquad\text{and}\qquad
	a(\mu')\subseteq \mathrm{int}(H^+)=\{ x\in \R^d : \langle x,w\rangle >\alpha \}. 
\end{equation}

\medskip
Let $E$ be the $d$-dimensional affine hyperplane in $\R^{d+1}$ defined by $E:=\{(x,1)\in\R^{d+1} : x\in\R^d\}$.
Denote by $e\colon \R^d\longrightarrow \R^{d+1}$ the embedding of $\R^d$ into $\R^{d+1}$ given by $x\longmapsto (x,1)$ for $x\in \R^d$. 
Now the affine map $a\colon \Delta_N\longrightarrow \R^d$ induces an affine map $e\circ a\colon \Delta_N\longrightarrow\R^{d+1}$.
The $e\circ a$ images of the vertices of $\Delta_N$ induce $1$-dimensional linear subspaces $L_{v_i}:=\spann  \{ (a(v_i),1)\}$ of $\R^{d+1}$ that pierce $E$ at the points $(e\circ a)(v_i)=(a(v_i),1)$, that is $L_{v_i}\cap E=\{(a(v_i),1)\}$, where $0\leq i\leq N$.

\medskip
Let $w':=(w,-\alpha)\in \R^{d+1}$.
Consider the $d$-dimensional affine hyperplane $E'$ in $\R^{d+1}$ defined by
\[
E':=\{ y\in \R^{d+1} : \langle y,w'\rangle =1 \}.
\]
Then the lines $L_{v_i}$ pierce $E'$ at the points
\begin{equation}
	\label{definition of b}
\widehat{a}(v_i):=\frac{1}{\langle a(v_i),w\rangle-\alpha}(e\circ a)(v_i)= \frac{1}{\langle a(v_i),w\rangle-\alpha}(a(v_i),1),
\end{equation}
for $0\leq i\leq N$.
Hence, we define a new affine map $\widehat{a}\colon \Delta_N\longrightarrow E'$ by its values on the vertices on $\Delta_N$. 
The choice we made in \eqref{separation} implies that the scaling factor $\frac{1}{\langle a(v_i),w\rangle-\alpha}$ in the definition of $\widehat{a}(v_i)$ is negative if and only if $v_i\in\vset(\mu)$, and positive if and only if $v_i\in\vset(\mu')$.

\medskip
Now apply to the affine map $\widehat{a}\colon \Delta_N\longrightarrow E'$ either the classical Tverberg theorem \cite[Thm.\,1]{Tverberg1966} or, with an appropriate coloring of the vertices of $\Delta_N$ and the assumption that $r$ is a prime, the Optimal colored Tverberg theorem \cite[Thm.\,2.1]{Blagojevic2009} \cite[Thm.\,2.1]{Blagojevic2011-2}.
We obtain $r$ pairwise disjoint non-empty faces, or rainbow faces (depending on the theorem we have used) $\sigma_1,\dots,\sigma_r$ of the simplex $\Delta_N$ such that $\widehat{a}(\sigma_1)\cap\dots\cap \widehat{a}(\sigma_r)\neq\emptyset$.
More precisely, there exist
\begin{compactitem}[\quad---]
\item pairwise disjoint non-empty subsets $I_1,\dots, I_r$ of the index set $I$ such that $\sigma_j:= \conv\{v_i : i\in I_j \}$ for all $1\leq j\leq r$, and 
\item non-negative scalars $\alpha_0,\dots,\alpha_N$ such that $\sum_{i\in I_1}\alpha_i=\dots=\sum_{i\in I_r}\alpha_i=1$,	
\end{compactitem}
with the property that
\begin{equation*}
	\sum_{i\in I_1}\alpha_i \widehat{a}(v_i)=\dots=\sum_{i\in I_r}\alpha_i \widehat{a}(v_i) \ \in \ \widehat{a}(\sigma_1)\cap\dots\cap \widehat{a}(\sigma_r) \ \subseteq E' .
\end{equation*}
Thus, using the definition of $\widehat{a}$ given in \eqref{definition of b}, we have that
\[
\sum_{i\in I_1}\frac{\alpha_i}{\langle a(v_i),w\rangle-\alpha}\,(a(v_i),1)=\dots=\sum_{i\in I_r} \frac{\alpha_i}{\langle a(v_i),w\rangle-\alpha}\,(a(v_i),1).
\]
Splitting the last equalities into the equalities of the first $d$ coordinates, and the equalities of the last $(d+1)$st coordinate we obtain that 
\begin{equation*}
	\label{first eq}
	\sum_{i\in I_1}\frac{\alpha_i}{\langle a(v_i),w\rangle-\alpha}\,a(v_i)=\dots=\sum_{i\in I_r} \frac{\alpha_i}{\langle a(v_i),w\rangle-\alpha}\,a(v_i) \ =: \ u \ \in \ \R^d,
\end{equation*}
and 
\begin{equation}
	\label{second eq}
	\sum_{i\in I_1}\frac{\alpha_i}{\langle a(v_i),w\rangle-\alpha}=\dots=\sum_{i\in I_r} \frac{\alpha_i}{\langle a(v_i),w\rangle-\alpha} \ =: \ \beta \ \in \ \R .
\end{equation}
The assumption that $\dim(\mu)\leq r-2$ yields the existence at lease one index set $I_{\ell}$ with the property $\vset(\mu)\cap\{v_i : i\in I_{\ell}\}=\emptyset$.
Therefore, the choice made in \eqref{separation} implies that $\beta=\sum_{i\in I_{\ell}}\frac{\alpha_i}{\langle a(v_i),w\rangle-\alpha} >0$, because all the summands are non-negative.
Consequently, from  \eqref{second eq} the point in $\aff(a(\sigma_1))\cap\dots\cap \aff(a(\sigma_r))$ we are looking for is
\[
\frac{1}{\beta}u= \sum_{i\in I_1}\frac{\alpha_i}{\beta(\langle a(v_i),w\rangle-\alpha)}\, a(v_i)=\dots=\sum_{i\in I_r} \frac{\alpha_i}{\beta(\langle a(v_i),w\rangle-\alpha)} \,a(v_i).
\]
Indeed, according to \eqref{separation} and the fact that $\beta>0$ holds we get that
\[
\big(v_i\in\vset(\mu)\  \Longrightarrow \ \frac{\alpha_i}{\beta(\langle a(v_i),w\rangle-\alpha)} \leq 0\big)
\quad\text{and}\quad 
\big(v_i\in\vset(\mu')\  \Longrightarrow \  \frac{\alpha_i}{\beta(\langle a(v_i),w\rangle-\alpha)} \geq 0\big).
\]
This completes the proof of both theorems, and gives a blueprint for proving all the remaining plus minus analogues directly from the corresponding Tverberg type theorems.

\providecommand{\bysame}{\leavevmode\hbox to3em{\hrulefill}\thinspace}
\providecommand{\href}[2]{#2}

\end{document}